\documentclass[12pt,reqno]{amsart}
\newtheorem{thm}{Theorem}[section]
\newtheorem{defi}{Definition}[section]

\newtheorem{cor}{Corollary}[section]
\newtheorem{pr}{Proposition}[section]
\theoremstyle{definition}
\newtheorem{rem}{Remark}[section]
\newtheorem{exm}{Example}[section]
\newcommand{\be}{\begin{equation}}
\newcommand{\ee}{\end{equation}}
\newcommand{\bea}{\begin{eqnarray}}
\newcommand{\eea}{\end{eqnarray}}
\newcommand{\beb}{\begin{eqnarray*}}
\newcommand{\eeb}{\end{eqnarray*}}
\usepackage{amssymb,url,amsfonts,amsthm,setspace,indentfirst,mathrsfs}
\usepackage{latexsym,amssymb,amsthm,amsmath}
\usepackage{graphicx}
\usepackage{minibox}
\usepackage{enumitem}
\usepackage{xcolor}
\usepackage[a4paper, margin=1in]{geometry} 
\usepackage{tikz}
\usepackage{amsmath}
\usepackage{amssymb}
\usetikzlibrary{positioning, arrows.meta, calc, decorations.markings}
\usepackage{mathtools}
\numberwithin{equation}{section}
\begin{document}
%
\title[On Lorentzian torse-forming structure manifolds]{ On Lorentzian torse-forming structure manifolds}
\author[A. A. Shaikh]{Absos Ali Shaikh}
\date{\today}
\address{\noindent\newline Department of Mathematics,
	\newline University of Burdwan, 
	\newline Golapbag, Burdwan-713104,
	\newline West Bengal, India} 
\email{aask2003@yahoo.co.in, aashaikh@math.buruniv.ac.in}

\dedicatory{}
\begin{abstract}
The main objective of the paper is to introduce a type of non-flat Lorentzian manifolds, called \emph {Lorentzian torse-forming structure manifolds (briefly, $LTFSMs$)}, which is a large class containing the Lorentzian concircular structure manifolds (briefly, $LCSMs$), LSP-Sasakian manifolds, and special type of LP-Sasakian manifolds with coefficient $\alpha$, LP-$\alpha$-Sasakian manifolds, LP-Sasakian manifolds and LP-$\alpha$-Kenmotsu manifolds. The interrelations among such various classes of manifolds are shown in Fig 1 and Fig 2.  Then we have investigated geometric properties of \emph{LTFSM} such as the expansion scalar (divergence of $\xi$) is non-zero and non-constant and the flow of $\xi$ is irrotational. In an $LCSM,$ the Ricci operator commutes with structure tensor $\phi$ \cite{Shaikh09}. But it is shown that in an \emph{LTFSM} the Ricci operator does not generally commute with the structure tensor $\phi$. However, if the gradient of the torse-forming scalar $\alpha$, i.e., $grad \ \alpha$ is collinear with the torse-forming vector field $\xi$, then its Ricci operator commutes with $\phi$, and it reduces to an \emph{LCSM}. We mention that a conformally flat $LCSM$ is quasi-Einstein hence is a manifold of quasi constant curvature. But we have  shown that a conformally flat \emph{LTFSM} is a generalized quasi-Einstein manifold in the sense of Chaki \cite{Chaki01} and also it is manifold of generalized quasi constant curvatue. Some curvature related properties of \emph{LTFSM} are deduced, and sectional curvature is determined. Then the existence of such a notion is ensured by several proper non-trivial examples. We have also explored the examples, which are both \emph{LCSM} and \emph{LTFSM}, only \emph{Lorentzian torse-forming manifold} (briefly, \emph{LTFM}), which characterizes twisted cubic spacetime, and also \emph{GRW} spacetime that characterizes \emph{LCSM}. Also, we have exhibited examples of \emph{LTFSM} in which the vector field $\xi$ is a unit timelike torse-forming but not concircular (see, Example $5.6-5.8$).

\end{abstract}
%
\maketitle
\section*{0. Introduction}

Let $M$ be an $n$-dimensional semi-Riemannian manifold endowed with a semi-Riemannian metric $g$ of signature $(p, n-p)$. If $p=0$ or $n$, then $M$ is called a Riemannian manifold. If $p=1$ or $n-1$, then $M$ is said to be a Lorentzian manifold. 
Thus the signature of the Lorentzian metric is either $(+, +, \dots, +, -)$ or $(-, +, \dots, +)$.
We note that the metric $g$ induces an inner product of the tangent space $T_pM$ at each point $p \in M$. Thus $g_p : T_pM \times T_pM \to \mathbb{R}$ is an inner product at each $p \in M$ such that $g_p(v, v)$ determines the length of the tangent vector $v \in T_pM$.
We also note that the length of a non-zero vector with respect to the 
Riemannian metric cannot be negative as the Riemannian metric assigns a 
positive definite inner product to the tangent space of the manifold. 
But the length of a zero vector with respect to the Riemannian metric 
is zero. Hence the length of a vector with respect to the Riemannian metric is always 
a non-negative real number. 
Also, the length of a vector at two different points of a Riemannian 
manifold can be same or different as the inner product is point dependent. 
But for an Euclidean metric the length of a vector is always same at all 
points of the space. 
Again, the length of a vector in a Lorentzian manifold can be positive, 
negative or zero. Even the length of a non-zero vector can be zero, which 
is called a null vector. If $g_p(v,v) > 0$ (resp., $< 0$), then $v$ is 
called a spacelike (resp., timelike). If $g_p(v,v) = 1$ or $-1$, 
then $v$ is called a unit spacelike or timelike vector. 
The geometry of the Lorentzian manifold is called Lorentzian geometry 
which plays a crucial role in the theory of general relativity and 
cosmology. Especially, the null cone exhibits an important role in 
this theory.\\
\indent Let $M$ be an $n$-dimensional semi-Riemannian manifold and $\nabla$ be the Levi-Civita connection of $g$. Then the Riemann, Ricci and scalar curvature are respectively defined as follows:
\begin{align*}
R(u,v)w &= \nabla_u \nabla_v w - \nabla_v \nabla_u w - \nabla_{[u,v]} w \\
S(u,w) &= \text{Trace } \{ v \mapsto R(v,u)w \} \quad \text{and} \\
r &= \text{Trace } S = \sum_{i=1}^n S(e_i, e_i),
\end{align*}
where $\{e_i : 1 \leq i \leq n\}$ is an orthonormal basis of the tangent space $T_pM$ at $p \in M$. Also we have $S(u,v) = g(Qu, v)$ where $Q$ is called the Ricci operator, which is an $(1,1)$ Ricci tensor and $S$ is a $(0,2)$ symmetric tensor. For $n=2$ the notion of $R,~ S$ and $r$ are equivalent in $M$ and hence throughout the paper we will consider $M$ of dimension $n \geq 3$. 
We note that the Ricci curvature in a certain direction is the average of the sectional curvatures of all planes containing that direction and it determines how the volume of a small ball of geodesics changes over time in comparison with an identical ball in flat space. In general relativity Ricci tensor describes the curvature of the spacetime caused by energy and momentum of the spacetime. A positive Ricci tensor causes geodesics to converge (decreasing volume), while a negative tensor causes them to diverge.
Again, the scalar curvature is a 
fundamental geometric measurement 
which describes how much the 
geometry of a curved space deviates 
from flat Euclidean space at a 
specific point. It provides a single 
``average" value that quantifies how 
the volume of a small ball 
in curved space differs from a 
sphere in flat space. In general relativity, scalar 
curvature represents the ``matter density" 
that sources gravity. It appears 
in the Einstein-Hilbert action where 
integrating the scalar curvature over 
a spacetime volume yields the total 
gravitational energy of the system.
The scalar curvature measures the 
volume of a small geodesic 
ball and its deviation from flat 
space depends directly on the 
scalar curvature. If $r > 0$, then the 
ball has less volume than a flat 
Euclidean ball of the same radius, 
and if $r < 0$, then it indicates a 
larger volume than a flat Euclidean 
ball of the same radius. In fact, 
$r$ is the exact arithmetic mean 
of all directional sectional 
curvatures in the manifold. In differential 
geometry, scalar curvature heavily
constraints the global shape of a manifold. For example, the Yamabe problem describes whether any smooth compact manifold can be conformally deformed to have a constant scalar curvature, which underpins many aspects of conformal field theory. \\
\indent A Lorentzian manifold is subject to extremely strict geometric restrictions when a global unit timelike torse-forming vector field is admitted. Such a structure inevitably forces the manifold to assume the form of a twisted cubic spacetime, realized as a twisted product of an one-dimensional interval and a Riemannian fiber, a result rigorously established by Mantica and Molinari \cite{MM17}. When the associated one-form is closed, or equivalently when the vector field is concircular, the geometry further specializes to that of a generalized Robertson--Walker spacetime, thereby furnishing a natural setting for cosmological models that possess a well-defined cosmic time and exhibit a smooth, homogeneous expansion of the spatial slices. From a purely topological perspective, the existence of a global continuous timelike vector field is by no means automatic; it requires the underlying manifold to be time orientable. In the compact case, this condition forces the Euler characteristic to vanish, while more generally, it guarantees a smooth foliation of the manifold into a sequence of spacelike hypersurfaces, ensuring a coherent global time function.

The kinematic properties of the flow defined by an unit timelike vector field are profoundly simplified under the torse-forming condition. An immediate consequence is the complete absence of acceleration, meaning that observers following the integral curves move along geodesic, freely falling trajectories without any external force. Furthermore, both the shear and rotation tensors vanish identically, indicating that the fluid elements undergo neither distortion nor rotational motion. Instead, the expansion scalar is found to be directly proportional to the associated structural function of the torse-forming field, reflecting a perfectly isotropic expansion in every spatial direction. This highly symmetric behavior encapsulates the essence of a homogeneous and isotropic cosmic fluid at the kinematic level.
Such a rigid kinematic structure imposes notable simplifications on the curvature tensors. The Ricci tensor reduces dramatically to a linear combination of the metric and the tensor product of the vector field with itself, with the coefficients determined entirely by the associated structural function and its derivatives along the flow. Consequently, the manifold necessarily acquires a quasi-Einstein character, wherein the curvature is governed solely by two scalar functions. When this curvature restriction is inserted into the Einstein field equations, the resulting total energy momentum tensor assumes the canonical form of a perfect fluid, characterized exclusively by an energy density and an isotropic pressure term. Thus, at the level of the total matter content, the geometry enforces a perfect fluid description.

However, this perfect fluid picture applies strictly to the aggregate matter sector. If one decomposes the total matter into distinct interacting components, such as a time-varying scalar field representing dark energy alongside conventional baryonic matter, the torse-forming framework provides a natural pathway to model more complex, imperfect fluid behavior. In this decomposed picture, the total energy-momentum tensor can be expressed as
\[
T_{\mu\nu} = T_{\mu\nu}^{\text{matter}} + T_{\mu\nu}^{\text{imperfect}}.
\]
The underlying structural constraints induce non-vanishing spatial gradients among these subcomponents. As a result, the uniform cosmic expansion gives rise to an effective bulk viscous pressure, which acts as an internal dissipative mechanism resisting the expansion. Moreover, whenever the gradient of the associated structural function fails to align perfectly with the observer field, a nontrivial heat flux emerges, thereby violating the perfect fluid condition for the individual matter constituents. This offers a rigorous geometric foundation for describing imperfect cosmological fluids, where dissipative effects such as bulk viscosity and energy transfer arise naturally from the geometry itself, rather than from adhoc phenomenological assumptions.
\vspace{1em}

In pseudo-Riemannian geometry and mathematical physics, the classification of spacetimes admitting highly structured timelike vector fields has developed along two largely independent directions. The first is primarily concerned with metric and topological properties, particularly warped and twisted product constructions such as generalized Robertson–Walker (GRW) spacetimes and twisted product manifolds. The second is centered on algebraic and geometric structures imposed on the manifold, including contact, paracontact, Sasakian, and concircular tensor frameworks. In some recent studies, these two perspectives have been treated as equivalent on the basis of common properties of their associated vector fields. Such an identification can obscure important distinctions between the underlying geometric categories. The present framework seeks to clarify this distinction and establish a more precise separation between metric classifications and additional structural geometries.\\
\indent The main objective of this work is to introduce and systematically study a broad category of connected and time oriented Lorentzian manifolds termed as \emph{Lorentzian torse-forming structure manifolds} ($LTFSMs$). This newly defined class is sufficiently general to encompass, as its particular subclasses, several well-known geometrical structures, including Lorentzian concircular structure manifolds, $LSP$-Sasakian manifolds and a special type of  $LP$-$\alpha$-Sasakian manifolds, $LP$-Sasakian manifolds with coefficient $\alpha$, $LP$-Sasakian manifolds, and  $LP$-Kenmotsu manifolds.  Following the introduction, we carry out a detailed investigation into the intrinsic geometric properties of $LTFSMs$ and establish the precise conditions under which such a structure reduces to the aforementioned various types of special manifolds.
The schematic diagram among these classes of manifolds are shown in Fig 1, where the implication notation means that such a class produces the other class of manifolds but not conversely. Also the inclusion relations among all such manifolds are exhibited in Fig 2. 

One of the key findings reveals that, in an $LTFSM$, the Ricci operator fails to commute with the structure tensor $\phi$. Nevertheless, when the gradient of the scalar function $\alpha$ is collinear with the torse-forming vector field, the Ricci operator does commute with $\phi$, and under this condition, the $LTFSM$ reduces precisely to an $LCSM$. Additionally, we demonstrate that any conformally flat $LTFSM$ satisfies the defining criteria of a generalized quasi-Einstein manifold, as introduced by Chaki \cite{Chaki01}. Also, we have shown that a conformally flat $LTFSM$ is a manifold of generalized quasi-constant curvature whereas an $LCSM$ is a manifold of quasi-constant curvature. Again, we have determined the associated scalars of an $LTFSM$ which is quasi-Einstein and generalized quasi-Einstein. Furthermore, several curvature properties of $LTFSMs$ are derived, including an explicit determination of the sectional curvature. To validate the theoretical framework, we construct several non-trivial and proper examples that ensure the consistency and existence of the proposed notion. In particular, the examples explored encompass a connected $LTFSM$ equipped with an $LCSM$ structure, as well as a standalone $LTFM$ that effectively characterizes a twisted cubic along with a $GRW$ spacetime. Moreover, we present a distinctive example of an $LTFSM$ wherein the vector field $\xi$ is a unit timelike torse-forming vector field which is not concircular.\\
\indent The paper is organized as follows:
In Section~1, we studied several fundamental geometric properties of torse-forming vector fields. Then we define $LTFSM$ along with other different types of manifolds. Section~2 is devoted to the study of conformally flat Lorentzian torse-forming structure manifolds. In Section~3, we investigate quasi-Einstein and generalized quasi-Einstein Lorentzian torse-forming structure manifolds. The sectional curvature of Lorentzian torse-forming structure manifolds is inspected in Section~4. In Section~5, we present several examples illustrating Lorentzian torse-forming structure manifolds and concicular structure manifolds. Finally, Section 6 and 7 are respectively deal with the conclusion of the whole work and future scope of research.
\begin{figure}[htbp]
    \centering
    
    \begin{tikzpicture}[
        box/.style={draw, rectangle, minimum width=2.5cm, minimum height=0.8cm, align=center, inner sep=2mm, font=\small},
            arrow/.style={-Implies, double, double equal sign distance, thick, shorten >=2pt, shorten <=2pt},
            biarrow/.style={Implies-Implies, double, double equal sign distance, thick, shorten >=2pt, shorten <=2pt}
        ]
            \node[box] (ltfm) {LTFM};
            \node[box, right=4.5cm of ltfm] (tc) {Twisted Cubic \\ manifolds};
            
            \node[box, below=1.5cm of ltfm] (ltfsm) {LTFSM};
            \node[box, below=1.5cm of tc] (grw) {GRW \\manifolds};
            
            \node[box, right=0.7cm of ltfsm] (lcsm) {LCSM};
        
            \node[box, left=1.5cm of ltfsm] (lpsca) {LP-Sasakian manifolds \\with coefficient $\alpha$};
            \node[box, above=1.3cm of lpsca] (lpk) {LP-Kenmotsu \\ manifolds};
            \node[box, below=1cm of lpsca] (lps) {LP-Sasakian \\manifolds};
        
            \node[box, below=1.5cm of ltfsm] (lpas) {LP-$\alpha$-Sasakian\\manifolds};
            \node[box, right=1.0 cm of lpas] (lsps) {LSP-Sasakian\\manifolds ($\alpha^2=1$)};
        
            \draw[biarrow] (ltfm) -- (tc);
            \draw[arrow] (ltfm) -- (ltfsm);
            \draw[arrow] (tc) -- (grw);
            
            \draw[arrow] (ltfsm) -- (lcsm); 
            \draw[arrow] (grw) -- (lcsm);
            
            \draw[arrow] (ltfsm) -- (lpsca);
            \draw[arrow] (lpsca) -- (lpk);
            \draw[arrow] (lpsca) -- (lps);
            \draw[arrow] (ltfsm) -- (lpas);
            \draw[arrow] (lpas) -- (lsps);
            \draw[arrow] (lpas.west) -- (lps.east);

    \end{tikzpicture}

    \vspace{0.5cm}
    {\large Fig 1.} Schematic diagram of different classes of manifolds
    \vspace{1.5cm}

    \begin{tikzpicture}[
        box/.style={draw, rectangle, minimum width=2.5cm, minimum height=0.8cm, align=center, inner sep=2mm, font=\small},
                setrel/.style={font=\Huge} 
            ]
                \node[box] (ltfsm2) {LTFSM};
                
                \node[setrel, left=0.15cm of ltfsm2] (sub2) {$\subset$};
                \node[box, left=0.15cm of sub2] (lpsca2) {LP-Sasakian\\with coefficient $\alpha$};
                
                \node[setrel, right=0.15cm of ltfsm2] (sub1) {$\subset$};
                \node[box, right=0.15cm of sub1] (ltfm2) {LTFM};
                
                \node[setrel, right=0.15cm of ltfm2] (equiv1) {$\equiv$};
                \node[box, right=0.15cm of equiv1] (tc2) {Twisted Cubic\\manifolds};
            
                \node[box, above=1.5cm of lpsca2] (lpk2) {LP-Kenmotsu\\manifolds};
                
                \node[box] at (ltfsm2 |- lpk2) (lpas2) {LP-$\alpha$-Sasakian\\ manifolds};
                
                \node[setrel, right=0.15cm of lpas2] (supset) {$\supset$};
                \node[box, right=0.15cm of supset] (lsps2) {LSP-Sasakian\\manifolds};
            
                \node[setrel] at ($(lpsca2)!0.5!(lpk2)$) (cap1) {$\cap$};
                \node[setrel] at ($(ltfsm2)!0.5!(lpas2)$) (cap2) {$\cap$};
            
                        \node[box, below=1.5cm of lpsca2] (lps2) {LP-Sasakian\\manifolds};
                        
                        \node[box, below=1.5cm of ltfsm2, minimum width=4.5cm] (lcsm2) {LCSM};
                        \node[box, minimum width=4.5cm] at (tc2 |- lcsm2) (grw2) {GRW\\manifolds};
                    
                        \node[setrel] at ($(lpsca2)!0.5!(lps2)$) (cup1) {$\cup$};
                        \node[setrel] at ($(ltfsm2)!0.5!(lcsm2)$) (cup2) {$\cup$};
                        \node[setrel] at ($(tc2)!0.5!(grw2)$) (cup3) {$\cup$};
                    
                        \node[setrel] at ($(lcsm2)!0.5!(grw2)$) (sub3) {$\subset$};

    \end{tikzpicture}

    \vspace{0.5cm}
    {\large Fig 2.} Inclusion diagram of different classes of manifolds

\end{figure}
\newpage
\section{Lorentzian torse-forming structure manifolds}
Let \(M\) be a connected Lorentzian manifold of dimension \(n\) equipped with a Lorentzian metric \(g\).
Let \(\chi(M)\) be the Lie algebra of all smooth vector fields on \(M\) and \(C^\infty(M)\) be the ring of smooth functions on \(M\). Then it is well known that \(\chi(M)\) is a \(C^\infty(M)\)-module. A vector field \(\xi \in \chi(M)\) is said to be spacelike (resp. timelike, null) if \(g(\xi, \xi) > 0\) (resp. \(g(\xi, \xi) < 0\), \(g(\xi, \xi) = 0\)).\\
 If \(\xi \) is a torse-forming vector field on \(M\) then there is an $1$-form \(\tau \) on \(M\) and a scalar function \(\alpha \) such that
\begin{equation}\label{TVF}
\nabla_X \xi = \alpha X + \tau(X)\xi,
\end{equation}
 for all \(X \in \chi(M)\), \(\tau(X) = g(X, \theta)\) and \(\nabla \) being the Levi-Civita connection of \(g\). The vector field \(\theta \) is called the associated vector field of the 1-form \(\tau \), and the scalar function $\alpha$ is said to be the associated structural function or torse-forming scalar of the vector field $\xi$.
In \eqref{TVF}, if we consider:\\
(i) \(\alpha = 1, \tau = 0\), then \(\xi \) is concurrent,\\
(ii) \(\tau = 0\) (or $\theta$ is collinear with $\xi$),  then $\xi$ is concircular,\\
(iii) \(\alpha = 0\), then \(\xi \) is recurrent.\\
(iv) \(\alpha = 0, \tau = 0\), then \(\xi \) is parallel.\\
\indent Let $M$ be a connected and time oriented Lorentzian manifold of dimension $n$, $n\geq3,$ endowed with a Lorentzian metric $g$ such that $M$ admits a torse-forming $\xi$. Then \eqref{TVF} holds. 
 \begin{equation}\label{TVF12}
     \nabla_X \xi = \alpha X + \tau(X) \xi
 \end{equation}
 If we apply a conformal change to the metric, $\bar{g} = e^{2\sigma} g$, where $\sigma \in C^\infty(M)$, the corresponding Levi-Civita connection $\bar{\nabla}$ is derived (using the Koszul formula) as follows.
 \begin{equation}
     \bar{\nabla}_X Y = \nabla_X Y + (X\sigma)Y + (Y\sigma)X - g(X,Y) \emph{grad} \sigma.
 \end{equation}
 
 Substituting $Y = \xi$, we obtain
 \begin{equation}
     \bar{\nabla}_X \xi = \nabla_X \xi + (X\sigma)\xi + (\xi\sigma)X - g(X,\xi) \emph{grad} \sigma,
 \end{equation}
 which yields by \eqref{TVF}
 \begin{align}
     \bar{\nabla}_X \xi &= (\alpha + (\xi\sigma))X + (\tau(X) + (X\sigma))\xi - g(X,\xi) \emph{grad} \sigma. \label{eq:conformal_xi}
 \end{align}
 This implies that torse-forming vector field is not invariant under a conformal transformation. However, if the conformal factor's gradient is pointwise collinear with $\xi$, then \eqref{eq:conformal_xi} entails that the torse-forming vector field $\xi$ is invariant under a conformal transformation. The same is true for a concircular vector field.\\
\indent If \(\tau \) is gradient, then there is a \(h \in C^\infty(M)\) such that \(\tau = dh\), and hence 
by Poincaré lemma \(d\tau = d^2h = 0\), and consequently \(\tau \) is a closed $1$-form.
 Again, if \(\tau \) is closed, then
\begin{align}
&~d\tau (X,Y)=0\notag\\
&\implies (\nabla _{X}\tau )(Y)=(\nabla _{Y}\tau )(X)\notag\\
&\implies g(\nabla _{X}\xi ,Y)=g(X,\nabla _{Y}\xi )\notag\\
&\implies g(\alpha X+\tau (X)\xi ,Y)
=g(X,\alpha Y+\tau (Y)\xi )
\quad \text{(by \eqref{TVF})}\notag\\
&\implies \tau (X)\eta (Y)=\tau (Y)\eta (X),\label{eq:taueta}\\
&\hspace{1.2em}\text{where } \eta (X)=g(X,\xi ),\quad \forall X\in\chi(M)\notag\\
&\implies \tau (X)=\frac{\tau (\xi )}{\|\xi \|^{2}}\eta (X)
=\sigma_1 \eta (X),\qquad
\sigma_1 =\frac{\tau (\xi )}{\|\xi \|^{2}},\quad
\forall X\in\chi(M).\notag
\end{align}
 i.e., 
\begin{equation}\label{Cod}
 \tau(X) = \sigma_1 \eta(X).
 \end{equation}
 
 This implies that the associated vector field $\theta$ of the $1$-form \(\tau \) and the associated vector field \(\xi \) of the $1$-form \(\eta \) are co-directional. It is entirely bound to the scalar function $\alpha$ and the metric projection of the field itself. Hence \eqref{TVF} takes the form\\
 \begin{equation}\label{TVF1}
 \nabla _{X}\xi =\alpha X+\sigma_1 \eta (X),
 \end{equation}
 where
 \begin{equation}\label{cof1}
 \sigma_1 = \frac{\tau(\xi)}{\Vert{}\xi\Vert{}^2}.
 \end{equation}
 If \(\xi \) is a torqued vector field, then \(\tau(\xi) = g(\xi, \theta) = 0\) and hence \eqref{cof1} yields \(\sigma_1 = 0\). Therefore \eqref{Cod} entails that \(\tau(X) = 0, \forall X\).
  Consequently,
 \begin{equation}\label{concir}
 \nabla _{X}\xi =\alpha X, \alpha \in C^{\infty}(M).
 \end{equation}
 Therefore \(\xi \) reduces to a concircular vector field.\\
 \indent  A specific subclass of torse-forming vector fields called, proper torqued vector fields can never be of constant length, because their lengths change across the manifold and they are fundamentally non unit fields.
For example in $\mathbb{R}^n$, the standard position vector field $X = (x_1, \dots, x_n)$ satisfies $\nabla_Y X = Y$. This is a specific type of torse-forming vector field (a concurrent vector field where $\alpha = 1$ and $\tau = 0$). Its length is $\sqrt{\sum_i x_i^2}$ which varies continuously and is explicitly not a unit vector field. Also non-trivial, non-unit torse-forming vector fields can be constructed on $m$-sphere $S^m$ by projecting a fixed Euclidean coordinate vector field onto the sphere's tangent bundle.\\ 
\indent Again, certain specialized subclasses such as anti-torqued vector fields are forced by their geometric constraints to always be unit geodesic vector fields. 
 Hence if \(\xi \) is a unit vector field, then \(g(\xi, \xi) = 1\) (for spacelike), or \(g(\xi, \xi) = -1\) (for timelike).\\
  Thus \eqref{TVF12} yields
  $$\tau (X)=\tau (\xi )\eta (X)\text{\ or\ }\tau (X)=-\tau (\xi )\eta (X)$$
  $$\text{i.e.,}~ \tau(X) = \beta \eta(X) \text{(for spacelike)}~ \text{or}~ \tau(X) = -\beta \eta(X), \text{(for timelike)},~ \beta = \tau(\xi).$$
So in this case, we have (from \eqref{TVF}), \(\nabla _{X}\xi =\alpha X+\beta \eta (X)\xi \) (for spacelike)
$$\text{or}~ \nabla _{X}\xi =\alpha X-\beta \eta (X)\xi \quad (\text{for\ timelike}).$$
Since \(\xi \) is a unit vector field, we have
 \(g(\xi, \xi) = \pm 1\) and hence \(g(\nabla_X \xi, \xi) = 0\),
  which yields, \(g(\alpha X + \beta \eta(X) \xi, \xi) = 0\) (for spacelike), or  \(g(\alpha X - \beta \eta(X) \xi, \xi) = 0\) (for timelike)
 i.e.,
 \begin{equation}\label{timelike}
\tau(X) = - \alpha \eta(X) \text{(for spacelike) or~} \tau(X) =  \alpha \eta(X) \text{(for timelike)}
 \end{equation}
  Hence if \(\xi \) is a unit timelike torse-forming vector field, then we get (by virtue of \eqref{timelike})
 \begin{equation}\label{UTTFV1}
 \nabla_X \xi = \alpha(X + \eta(X) \xi),
 \end{equation}\label{timelike2}
  and if \(\xi \) is a unit spacelike torse-forming vector field, then
 \begin{equation}\label{UTTF}
 \nabla_X \xi = \alpha(X - \eta(X) \xi).
 \end{equation}
 
 Hence we note that if $\xi$ is a unit timeline (resp. spacelike) torse-forming vector field such that $grad~ \alpha$ is collinear with $\xi$, then $\xi$ turns out to be a proper concircular vector field.
  Therefore if a Lorentzian manifold \(M\) admits a unit timelike torse-forming vector field \(\xi \), then \eqref{UTTFV1} holds, which
   implies that
   \begin{equation}\label{UTT1}
  (\nabla_X \eta)(Y) = g(\nabla_X \xi, Y),
  \end{equation}
 i.e.,
 \begin{equation}\label{UTT2}
  (\nabla_X \eta)(Y) = \alpha [g(X, Y) + \eta(X) \eta(Y)],
 \end{equation}
   and hence
    $$(\nabla_X \eta)(Y) = (\nabla_Y \eta)(X)~ \text{for all}~ X, Y \in \chi(M).$$
 This implies $d\eta(X, Y) = 0$ and so $\eta$ is a closed $1$-form.
  From \eqref{UTT1} and \eqref{UTT2} we get
\begin{equation}\label{div}
 div \text{ }\xi = \alpha(n-1) \neq 0 \neq \text{constant}
\end{equation}
   and $g(\nabla_X \xi, Y) = g(\nabla_Y \xi, X) \implies \text{curl } \xi = 0.$
   Thus \(\xi \) is an irrotational vector field.
 If \(\xi = grad~  f\), \(f \in C^{\infty}(M)\), then \(\exists \) a vector field \(V\) on \(M\) such that \cite{MC10}
 \begin{equation}\label{grad1}
 \xi = e^f V
  \implies \nabla_X V = \beta X,
 \end{equation}
 where $\beta = \alpha e^f, \beta \in C^{\infty}(M).$ 
 As \(\eta(X) = g(X, \xi) = g(grad~ f, X) = (Xf)\).
  Hence
 \begin{equation}\label{concir1}
 \nabla_X \xi = \alpha(X + \eta(X)\xi) = \alpha(X + (Xf)\xi).
 \end{equation}
  This implies that \(V\) is a concircular vector field. If \(V\) is concircular, then every vector field collinear with \(V\) is a torse-forming vector field.
From \eqref{UTTFV1}, it follows that \(\xi \) is a unit timelike concircular vector field if it also satisfies
\begin{equation}\label{concir2}
\nabla_X \alpha = (X \alpha) = \mu \eta(X), \text{ } \mu \in C^{\infty}(M).
\end{equation}
In view of the above discussions, we can state the following:
\begin{pr}
Let \(M\) be a Lorentzian manifold of dimension \(n, n \geqslant 3\), such that \(\xi \) is a unit timelike torse-forming vector field on \(M\). Then the following statements hold:\\
(i) If \(\xi \) is a torqued vector field and \(\tau \) is a gradient form, then \(\xi \) is a concircular vector field.\\
(ii) If \(\xi \) is gradient, then \(\exists \) a concircular vector field \(V\) such that every vector field collinear with \(V\) is a torse-forming vector field.\\
(iii) The relations \eqref{UTTFV1} and \eqref{UTT2} hold.\\
(iv) $div$ \( \xi\) is non-zero and non-constant. If \(\xi \) is recurrent, then $div$ \( \xi = 0\).\\
(v) \(\xi \) is a irrotational vector field as $curl$ \( \xi = 0\).\\
(vi) If \eqref{concir2} holds, then \(\xi \) is a unit timelike concircular vector field.
\end{pr}
\indent If we set
\begin{equation}\label{Omega}
\Omega(X, Y) = g(\phi X, Y) = \frac{1}{\alpha}(\nabla_X \eta)(Y)
\end{equation}
 for all \(X, Y \in \chi(M)\), then \eqref{UTT2} entails that
 \begin{equation}\label{timelike1}
 \phi X = X + \eta(X) \xi.
 \end{equation}
This implies that a unit timelike torse-forming vector field on \(M\) produces an (1,1) tensor field \(\phi \).
From \eqref{timelike1} we get
\begin{equation}\label{TFV1}
g(\phi X, \phi Y) = g(X, Y) + \eta(X)\eta(Y)
\end{equation}
\begin{equation}\label{TFV2}
\phi^2 X = X + \eta(X) \xi,~\phi \xi =0,~
 \Omega (X,Y)=g(\phi X,Y)=g(X,\phi Y)=\Omega (Y,X)
\end{equation}
and hence \(\phi \) is an (1,1) symmetric tensor field and \(\Omega \) is a symmetric $(0,2)$ tensor.
 This implies that a Lorentzian manifold admitting a unit timelike torse-forming vector field satisfying \eqref{Omega} is a special type of Lorentzian almost paracontact manifolds \cite{Sato77}. 
We note that the eigen values of $\phi$ are 0 of multiplicity 1, and 1 of multiplicity $(n-1)$ corresponding to the $(n-1)$-dimensional distribution orthogonal to $\xi.$
Hence there exists a local orthonormal frame $\{e_1,e_2,..e_{n-1},\xi\}$ at each point of the manifold, where $\{e_1,e_2,..e_{n-1},\xi\}$ are spacelike vectors orthogonal to $\xi$. i.e., $g(e_i,e_j)=\delta_{ij},~i,j\in\{1,2,...n-1\}$, $g(e_i,\xi)=0, i\in\{1,2,..n-1\}$, $g(\xi,\xi)=-1.$ In this basis, $\phi$ acts as a projection operator onto the $(n-1)$-dimensional distribution orthogonal to $\xi.$ And hence trace of $\phi = n-1 .$ Also $\phi(e_i)=e_i$ for spacelike vectors $e_i$, and  $\phi(\xi)=0$ for timelike vector $\xi$. Hence the matrix representation of $\phi$ is given as 
\[
{\phi} = \begin{pmatrix} 
I_{(n-1)\times(n-1)} & \dots & 0 \\ 
\vdots & \ddots & 0 \\ 
0 & \dots & 0 
\end{pmatrix}.
\]
Since $\phi X = X + \eta(X) \xi$, it is the orthogonal projection operator onto the spatial hypersurface orthogonal to the temporal flow line of $\xi.$ Hence $\phi$ determines the spatial geometry of the manifold, and it annihilates the time dimension isolating the local rest frame. Also $\phi$ acts as an identity map on any pure spatial vector preserving spatial lengths and angles. Since $\xi$ is torse-forming and $\phi$ is symmetric, its dual form $\eta$ is closed. Hence by Frobenious theorem, it follows that the distribution annihilated by $\eta$ (the range of $\phi$) is completely integrable. Geometrically, the manifold is sliced up cleanly like a deck of cards into smooth, nested spatial hypersurfaces (leaves). In general relativity, $\phi$ is used to breakdown 4-dimensional spacetime tensors into measurable 3-dimensional physical quantities. For example, energy-momentum tensor is decomposed using $\phi$ to separate the physical energy density, isotropic pressure, shear stress observed by a moving fluid element.    

From the previous discussion, we  now define various classes of manifolds as follows:
\begin{defi}
A generalized Robertson-Walker (GRW) manifold is defined as an $n$-dimensional ($n \ge 3$) connected time oriented Lorentzian manifold $M$ which is a warped product of the form:
\[
M = -I \times_{f} M^*
\]
where
     $I$ is an open interval of the real line $\mathbb{R}$ equipped with the negative definite metric $-dt^2$ and
     $M^*$ is an $(n-1)$-dimensional Riemannian manifold,
     $f$ is a smooth, strictly positive function of time on $I$, known as the warping function.
The metric tensor $g$ of $M$ is given locally by:
\[
ds^2 = -dt^2 + f(t)^2 g^*,
\]
where $g^*$ is the Riemannian metric of $M^*$. A $GRW$  spacetime is a connected time oriented $4$-dimensional manifold endowed with the metric 
 \[
 ds^{2} = -dt^{2} + f^{2}(t) \left( dx^{2} + dy^{2} + dz^{2} \right).
 \]
\end{defi}

 \begin{defi}\cite{MM17}
 A twisted cubic manifold $M$ is defined as an $n$-dimensional ($n \ge 3$) connected time oriented Lorentzian manifold which is a warped product of the form:
 \[
 M = -I \times_{f} M^*
 \]
 where
      $I$ is an open interval of the real line $\mathbb{R}$ equipped with the negative definite metric $-dt^2$ and
      $M^*$ is an $(n-1)$-dimensional Riemannian manifold,
      $f$ is a smooth, strictly positive function of time and position, known as the warping function.
 The metric tensor $g$ of $M$ is given locally by:
 \[
 ds^2 = -dt^2 + f(t,\mathbf{x})^2 g^*,
 \]
 where $g^*$ is the Riemannian metric of $M^*$.
 \end{defi}
 A twisted cubic spacetime is a connected time oriented $4$-dimensional manifold endowed with the metric 
 \[
 ds^{2} = -dt^{2} + f^{2}(t,\mathbf{x}) \left( dx^{2} + dy^{2} + dz^{2} \right).
 \]
 A $GRW$ spacetime is a rigid, specialized subclass of a more flexible twisted cubic spacetime.
\begin{defi}
A connected time oriented Lorentzian manifold \(M\) of dimension \(n, n \geq 3\), admitting a unit timelike torse-forming vector field is called a Lorentzian torse-forming manifold (briefly, LTFM).
\end{defi}
We note that Mantica and Molinari \cite{MM17} showed that an $LTFM$ is a twisted cubic manifold and vice-versa.
\begin{defi}\label{DLTFSM}
A connected time oriented Lorentzian manifold $M$ of dimension \(n, n \geq 3\), admitting a unit timelike torse-forming vector field satisfying the condition \eqref{Omega} is called a Lorentzian torse-forming structure manifold (briefly, LTFSM). The torse-forming vector field $\xi$ is said to be characteristic or flow vector field of the manifold, and the $(1,1)$ tensor field $\phi$ will be designated as structure tensor field of the manifold. We note that LTFSMs is a subclass of LTFMs as well as twisted cubic manifolds.
\end{defi}
\begin{defi}\label{DLCSM}
A connected time oriented Lorentzian manifold $M$ of dimension \(n, n \geq 3\), admitting a unit timelike concicular vector field is said to be a Lorentzian concircular structure manifold (briefly, $LCSM$) if it satisfies the relations \eqref{UTTFV1}, \eqref{UTT2}, \eqref{concir2} and \eqref{Omega}.
\end{defi}
We note that in an $LCSM$, the relations \eqref{timelike1}-\eqref{TFV2} hold. We also note that $LCSMs$ is a subclass of $GRW$ manifolds and hence $(CS)_4$  spacetimes is a subclass of $GRW$ spacetimes \cite{MM19}.
\begin{defi}\cite{DSS2002}\label{DLPSM}
A connected time oriented Lorentzian manifold \(M\), of dimension \(n, n \geq 3\), admitting a unit timelike torse-forming vector field such that \eqref{Omega} and \eqref{LPS} (given below) hold, is called an LP-Sasakian manifold with a coefficient \(\alpha \). In particular, if $\alpha=1,$ then it is an $LP$-Sasakian manifold.\end{defi}
\begin{equation}\label{LPS}
\begin{rcases}
\phi^2 X = X + \eta(X)\xi,~g(\phi X, \phi Y) = g(X, Y) + \eta(X)\eta(Y),
g(X,\xi)=\eta(X),\\ 
\Omega (X,Y)=g(\phi X,Y)=g(X,\phi Y)=\Omega (Y,X), \phi \xi=0, \eta(\phi X)=0, \eta(\xi) = -1,\\
\forall X,Y\in \chi(M).
\end{rcases}
\end{equation}
\begin{defi}
A Lorentzian para-$\alpha$-Sasakian (briefly, LP-$\alpha$-Sasakian) manifold $M$ of dimension \(n,~ n \geq 3\), is a connected manifold equipped with a Lorentzian almost paracontact metric structure $(\phi, \xi, \eta, g)$ satisfying relation \eqref{LPS} and for Levi-Civita connection $\nabla$ with respect to the Lorentzian metric $g$ realizes:
\begin{equation}\label{LPAS}
(\nabla_X \phi)Y = \alpha \big( g(X, Y)\xi + \eta(Y)X \big), 
\end{equation}
for all vector fields X, Y on M,
where $\alpha$ is a non-zero smooth function. If $\alpha = 1$, the structure reduces to an LP-Sasakian manifold.
\end{defi}
\begin{defi}\label{DLPKM} 
A Lorentzian para-$\alpha$-Kenmotsu (briefly, LP-$\alpha$-Kenmotsu) manifold of dimension \(n, n \geq 3\), is a connected manifold equipped with a Lorentzian almost paracontact metric structure $(\phi, \xi, \eta, g)$ satisfying the relation \eqref{LPS} and the Levi-Civita connection $\nabla$ with respect to the Lorentzian metric $g$ satisfies:
\begin{equation}\label{LPAK1}
(\nabla_X \phi)Y = \alpha(g(\phi X, Y)\xi - \eta(Y)\phi X).
\end{equation}
Especially, if $\alpha=1,$ then it becomes an LP-Kenmotsu manifold.
\end{defi}
\begin{defi}\label{DLSPSM}
A Lorentzian special para-Sasakian (briefly, LSP-Sasakian) manifold $M$ of dimension \(n, n \geq 3\), is a connected manifold equipped with a Lorentzian almost paracontact metric structure $(\phi, \xi, \eta, g)$ satisfying relations \eqref{UTT2} for $\alpha^2=1$, and \eqref{LPS}. 
\end{defi}
\begin{rem}
We note that in Definitions 1.6 - 1.8, there is no expression of $\phi$, only we have the expression of $\phi^2$. But in Definition 1.4, 1.5 and 1.9 we have both the same expressions of $\phi$ and $\phi^2$.  
\end{rem}
\begin{rem}
We also note that in Definitions 1.1 - 1.3, there is a concept of the $(1,1)$ tensor $\phi$ associated with the unit timelike torse-forming vector field $\xi$.
\end{rem}
Throughout the paper we will cosider the vector field $X,Y, Z,U... \in \chi(M),$\\ and $a,b,c, \alpha,\mu...\in C^{\infty}(M).$ We consider a manifold $M$ of dimension \(n, n \geq 3\), with the characteristic vector field $\xi$ as the unit timelike torse-forming. Hence in an $LTFSM$ the relations \eqref{UTT1}, \eqref{UTT2}, and \eqref{Omega}-\eqref{TFV2} hold.\\
\indent From \eqref{UTTFV1} we get
\begin{equation}\label{LCS}
\nabla_X \nabla_Y \xi = (X\alpha)[Y + \eta(Y)\xi] + \alpha \left[ \nabla_X Y + \nabla_X \eta(Y) \xi + \eta(Y) \nabla_X \xi \right]
\end{equation}
and 
\begin{equation}\label{LCS1}
\nabla_{[X,Y]} \xi = \alpha \left( [X,Y] - \eta([X,Y])\xi \right),
\end{equation}
where 
$$(X\alpha)=X(\alpha)=(\nabla_X\alpha)=d\alpha(X)=g(\text{$grad$ } \alpha, X).$$\\
In view of \eqref{LCS} and \eqref{LCS1}, we obtain (by virtue of \eqref{UTT2}) the following:\\
\begin{equation}\label{RC1}
R(X,Y)\xi =(X\alpha )[Y+\eta (Y)\xi ]-(Y\alpha )[X+\eta (X)\xi ]+\alpha ^{2}[\eta (Y)X-\eta (X)Y],
\end{equation}
\begin{equation}\label{eta}
\begin{split}
\eta (R(X,Y)Z) &= (Y\alpha )[g(X,Z)+\eta (X)\eta (Z)] - (X\alpha )[g(Y,Z)+\eta (Y)\eta (Z)] \\
&\quad + \alpha^{2}[g(Y,Z)\eta (X)-g(X,Z)\eta (Y)].
\end{split}
\end{equation}
From \eqref{RC1}, we get\\
\begin{equation}\label{Ricci}
S(Y,\xi )=-(n-2)(Y\alpha )+[(\xi\alpha )+\alpha ^{2}(n-1)]\eta (Y),
\end{equation}
\begin{equation}\label{Qxi}
Q\xi =-(n-2)\text{$grad$\ }\alpha +[(\xi \alpha )+\alpha ^{2}(n-1)]\xi, 
\end{equation}
where $S(X, Y) = g(QX, Y),$
 \(Q\) being the Ricci operator.
From \eqref{RC1}, we can obtain by virtue of \eqref{timelike1} that\\
\begin{equation}
R(\xi,\phi Z,X,Y)- R(\xi ,Z,X,Y)=0.
\end{equation}
Since $\phi \xi=0,$ by using \eqref{Ricci}, we get 
\begin{equation}\label{Ricci5}
S(\phi X,Y)-S(X,\phi Y)=(n-2)[(X\alpha)\eta(Y)-(Y\alpha)\eta(X)]
\end{equation}
and hence
\begin{equation}\label{Qphi}
Q(\phi X)-\phi (QX)=(n-2)[(X\alpha )\xi -\eta (X)\text{$grad$\ }\alpha ].
\end{equation}
If \(Q(\phi X) - \phi(QX) = 0\), then for \(n \geq 3\), we have from \eqref{Qphi} that
$$(X\alpha )\eta (Y)-\eta (X)(Y\alpha )=0,$$
which yields \((X\alpha )=-(\xi \alpha )\eta (X)=\rho\eta(X), \rho=-(\xi\alpha)\) 
and hence \(grad~\alpha =\rho \xi \).
Consequently, \(Q(\phi X) - \phi(QX) = 0\).
Thus we can state the following:
\begin{thm}
In an LTFSM of dimension $n$, \(n \geq 3\), the Ricci-operator $Q$ commutes with the structure tensor \(\phi \) if and only if the gradient of the torse-forming scalar function \(\alpha \) is collinear with the flow vector field \(\xi \). Hence \(\xi \) is a unit timelike concircular vector field and consequently the manifold reduces to a Lorentzian concircular structure manifold.
\end{thm}
\indent Physically, the torse-forming associated function $\alpha$ is exactly the Hubble expansion parameter $H(t)$ of the universe. The equation $\nabla_X \xi=\alpha(\phi X)$ means that any two neighbouring galaxies are moving apart at a speed proportional to the spatial distance, which is the exact definition of Hubble's law in an FLRW universe. Geometrically, $\phi$ is the spatial compass that isolates $3$-dimensional space from $4$-dimensional spacetime. Physically, $\phi$ forces the universe to expand uniformly and isotropically in all spatial directions without twisting or distorting, mimicking the exact large scale behaviour of our actual cosmos.  

\section{Conformally flat Lorentzian Torse-forming Structure Manifolds}
\indent We suppose that an $LTFSM$ of dimension \(n > 3\) is conformally flat. Then vanishing of conformal curvature tensor yields\\
\begin{align}\label{R04}
R(X,Y,Z,U)
=&\; \frac{1}{n-2}\Big[
S(Y,Z)g(X,U)-S(X,Z)g(Y,U) \notag\\
&\qquad +S(X,U)g(Y,Z)-S(Y,U)g(X,Z)
\Big] \notag\\
&-\frac{r}{(n-1)(n-2)}
\Big[
g(Y,Z)g(X,U)-g(X,Z)g(Y,U)
\Big],
\end{align}
for all vector fields \(X, Y, Z, U \in \chi(M)\).\\
In view of \eqref{eta}, \eqref{Ricci} and \eqref{R04} (after straightforward calculation), we obtain\\
\begin{equation}\label{Ricci1}
S(Y,Z)=\psi g(Y,Z)+\delta_1 \eta (Y)\eta (Z)+\nu_1[\zeta (Y)\eta (Z)+\zeta (Z)\eta (Y)],
\end{equation}
where \(\zeta(Y) = g(\text{$grad$ } \alpha, Y)=(Y\alpha)\) for all \(Y\) and 
$\psi, \delta_1$ and $\nu_1$ are scalar functions given by\\
$\psi=\frac{r}{n-1}-(\xi \alpha)-\alpha^2,~ \delta_1=\frac{r}{n-1}+(\xi \alpha)(n-4)-n\alpha^2,~\nu_1=(n-2).$\\
\indent We now project $grad~\alpha$ onto the spatial distribution orthogonal to $\xi$ as\\ $grad^\perp \alpha=grad~\alpha+(\xi\alpha)\xi$. Then $g(grad^\perp \alpha,\xi)=0$ and $||grad^\perp \alpha||^2=||grad~ \alpha||^2+(\xi\alpha)^2.$ Hence \eqref{Ricci1} takes the form
\begin{equation}\label{Ricci6}
S(Y,Z)=\psi g(Y,Z)+\delta \eta (Y)\eta (Z)+\nu[B (Y)\eta (Z)+B (Z)\eta (Y)],
\end{equation}
where 
$\psi, \delta$ and $\nu$ are scalar functions given by\\
$\psi=\frac{r}{n-1}-(\xi \alpha)-\alpha^2,~ \delta=\frac{r}{n-1}+n(\alpha^2+(\xi \alpha)),~\nu=\nu_1||grad^\perp \alpha||, B(Y)=g(\mu_2,Y)$ for all $Y,$ $\mu_2=\frac{grad^\perp \alpha}{||grad^\perp \alpha||}$ such that $||\mu_2||=1.$
If $\alpha$ is chosen as spatial proper distance slicing coordinate system, the metric functions are normalized such that the spatial gradient $grad^\perp\alpha$ has a norm 1 to represent an isometric embedding or a normalized direction of spatial inhomogeneity, and hence the universe cannot be perfectly isotropic. Consequently the spacetime must be twisted and the spatial gradient $grad^\perp\alpha$ acts as a directional vector pointing towards the spatial variations of the universe's expansion rate. Hence the spatial coordinates are explicitly calibrated such that travelling of 1 unit physical distance corresponds exactly to 1 unit change in the expansion parameter. Therefore $grad^\perp\alpha$ is a unit vector if and only if the spatial variations of the torse-forming scalar $\alpha$ across the hyper surfaces orthogonal to $\xi$ are perfectly normalized to a magnitude of 1.
\begin{defi}\cite{Chaki01}
A non-flat semi-Riemannian manifold $(M^n, g)$ of dimension $n \geq 3$ is called a \textit{generalized quasi-Einstein manifold}, if its Ricci tensor $S$ of type $(0, 2)$ is not identically zero and satisfies the relation:
\begin{equation}\label{Ricci4}
S(X, Y) = a g(X, Y) + b A(X)A(Y) + c [A(X)B(Y) + A(Y)B(X)],
\end{equation}
for any vector field $X$ and $Y$ on $M$, where $A$ and $B$ are nonwhere vanishing $1$-forms such that $A(X) = g(X, \mu_1)$ and $B(X) = g(X, \mu_2)$, $\mu_1$ and $\mu_2$ are mutually orthogonal unit vector fields, and $a, b, c\in C^{\infty}(M).$
\end{defi}
If we consider $\mu_1=\xi$ and $\mu_2=\frac{grad^\perp \alpha}{||grad^\perp \alpha||},$ then  from \eqref{Ricci6}, we can state the following:
\begin{thm}
A conformally flat LTFSM of dimension \(n > 3\) is a generalized quasi-Einstein manifold in the sense of Chaki.
\end{thm}
\begin{cor}
Any 3-dimensional $LTFSM$ is a generalized quasi-Einstein manifold in the sense of Chaki.
\end{cor}
\begin{cor}
A conformally flat $LCSM$ of dimension \(n > 3\) is a quasi-Einstein manifold.
\end{cor}
Again using \eqref{Ricci1} in \eqref{R04}, we obtain
\begin{equation}
\begin{aligned}
R(X,Y,Z,U)
=&\; a_1\,[g(Y,Z)g(X,U)-g(X,Z)g(Y,U)] \\
&+ \tilde {b_1}\,[g(Y,Z)\eta(X)\eta(U)+g(X,U)\eta(Y)\eta(Z) \\
&\qquad -g(Y,U)\eta(X)\eta(Z)-g(X,Z)\eta(Y)\eta(U)] \\
&+ \,[g(Y,Z)\{\zeta(X)\eta(U)+\zeta(U)\eta(X)\} \\
&\qquad +g(X,U)\{\zeta(Y)\eta(Z)+\zeta(Z)\eta(Y)\} \\
&\qquad -g(X,Z)\{\zeta(Y)\eta(U)+\zeta(U)\eta(Y)\} \\
&\qquad -g(Y,U)\{\zeta(X)\eta(Z)+\zeta(Z)\eta(X)\}] ,
\end{aligned}
\label{Rc04}
\end{equation}
where \(a_1=\frac{1}{n-2}[\frac{r}{(n-1)}-2(\alpha^2+(\xi \alpha))]\),
\(\tilde {b_1}=\frac{1}{n-2}[\frac{r}{n-1}+(n-4)(\xi \alpha)-n\alpha^2]\). Hence as before, \eqref{Rc04} can be rewritten as  
\begin{equation}
\begin{aligned}
R(X,Y,Z,U)
=&\; a_1\,[g(Y,Z)g(X,U)-g(X,Z)g(Y,U)] \\
&+  b_1\,[g(Y,Z)\eta(X)\eta(U)+g(X,U)\eta(Y)\eta(Z) \\
&\qquad -g(Y,U)\eta(X)\eta(Z)-g(X,Z)\eta(Y)\eta(U)] \\
&+ \,c_1[g(Y,Z)\{B(X)\eta(U)+B(U)\eta(X)\} \\
&\qquad +g(X,U)\{B(Y)\eta(Z)+B(Z)\eta(Y)\} \\
&\qquad -g(X,Z)\{B(Y)\eta(U)+B(U)\eta(Y)\} \\
&\qquad -g(Y,U)\{B(X)\eta(Z)+B(Z)\eta(X)\}] ,
\end{aligned}
\label{Rc05}
\end{equation}
where \(a_1=\frac{1}{n-2}[\frac{r}{(n-1)}-2(\alpha^2+(\xi \alpha))]\),
\(b_1=\frac{1}{n-2}[\frac{r}{n-1}+(n-6)(\xi \alpha)-n\alpha^2],~c_1=||grad^\perp \alpha||\).
\begin{defi}\cite{Chaki01}
A non-flat semi-Riemannian manifold \((M^n, g)\) \((n \geq 3)\) is called a manifold of generalized quasi-constant curvature if its curvature tensor \(R\) of type \((0,4)\) is of the form

\begin{align}\label{QCN}
R(X,Y,Z,W) &= a_2 \big[ g(Y,Z)g(X,W) - g(X,Z)g(Y,W) \big] \nonumber \\
&\quad + b_2 \big[ g(Y,Z)A(X)A(W) - g(X,Z)A(Y)A(W) \nonumber \\
&\quad \quad + g(X,W)A(Y)A(Z) - g(Y,W)A(X)A(Z) \big] \nonumber \\
&\quad + c_2 \big[ g(Y,Z)\{ A(X)B(W) + A(W)B(X) \} \nonumber \\
&\quad \quad - g(X,Z)\{ A(Y)B(W) + A(W)B(Y) \} \nonumber \\
&\quad \quad + g(X,W)\{ A(Y)B(Z) + A(Z)B(Y) \} \nonumber \\
&\quad \quad - g(Y,W)\{ A(X)B(Z) + A(Z)B(X) \} \big],
\end{align}
$X,Y,Z,W\in \chi(M),$ where $A$ and $B$ are nonwhere vanishing $1$-forms such that $A(X) = g(X, \mu_1)$ and $B(X) = g(X, \mu_2)$, $\mu_1$ and $\mu_2$ are mutually orthogonal unit vector fields, and $a_2, b_2, c_2\in C^{\infty}(M).$ In particular, if $c_2$=0, then \eqref{QCN} takes the form of a manifold of quasi-constant curvature. Also if $b_2=c_2=0$, then \eqref{QCN} reduces to the form of a manifod of constant curvature.\\
If we take $\mu_1=\xi$ and $\mu_2=\frac{grad^\perp \alpha}{||grad^\perp \alpha||}$, then  from \eqref{Rc05}, we can state the following:
\end{defi}
\begin{thm}
A conformally flat $LTFSM$ of dimension \(n > 3\) is a manifold of generalized quasi-constant curvature.
\end{thm}
\begin{cor}
Any 3-dimensional $LTFSM$ is a manifold of generalized quasi-constant curvature.
\end{cor}
 We note that if \eqref{concir2} holds then an $LTFSM$ reduces to an $LCSM$. Hence we can state the following:
 \begin{cor}
 A conformally flat $LCSM$ of dimension \(n > 3\) is a manifold of quasi-constant curvature.
 \end{cor}
 \begin{cor}\cite{DSJ}
  A conformally flat LP-Sasakian manifold of dimension \(n > 3\) with a non-constant coefficient \(\alpha \) is a manifold of quasi-constant curvature.
  \end{cor}
  \begin{cor}\cite{DSJ}
   A conformally flat LP-Sasakian manifold of dimension \(n > 3\) with a constant coefficient \(\alpha \) is a manifold of constant curvature.
   \end{cor}
   \begin{cor}
     A conformally flat LP-Sasakian manifold of dimension \(n > 3\) is of constant curvature.
    \end{cor}
   \begin{cor}
   A conformally flat LP-$\alpha$-Kenmotsu manifold of dimension \(n > 3\) with a non-constant coefficient \(\alpha \)  is a manifold of quasi-constant curvature.
   \end{cor}
   \begin{cor}
      A conformally flat LP-$\alpha$-Kenmotsu manifold of dimension \(n > 3\) with a constant coefficient \(\alpha \)  is a manifold of constant curvature.
      \end{cor} 
   \begin{cor}
      A conformally flat LP-$\alpha$-Sasakian manifold of dimension \(n > 3\) with a non-constant coefficient \(\alpha \) is a manifold of quasi-constant curvature.
      \end{cor}
     \begin{cor}
           A conformally flat LP-$\alpha$-Sasakian manifold of dimension \(n > 3\) with a constant coefficient \(\alpha \) is a manifold of constant curvature.
           \end{cor} 
\section{Quasi-Einstein and Generalized quasi-Einstein Lorentzian Torse-forming Structure Manifolds}
We first consider an $LTFSM$ of dimension \(n > 2\) which is quasi-Einstein. Then its Ricci tensor \(S\) of type \((0,2)\) is of the form\\
\begin{equation}\label{Ricci2}
S(X,Y) = a_3 g(X,Y) + b_3 \eta(X) \eta(Y)
\end{equation}
for all \(X, Y \in \chi(M)\), \(a_3, b_3 \in C^\infty(M)\), \(\eta(X) = g(X, \xi)\). From (3.1) we have
\begin{equation}\label{SC}
r = n a_3 - b_3.
\end{equation}
Again, from \eqref{Ricci2} we have
\begin{equation*}
S(\xi,\xi)=-a_3+b_3,
\end{equation*}
which yields by virtue of \eqref{Ricci} that 
\begin{equation}\label{Ricci6}
a_3-b_3=(n-1)[(\xi \alpha)+\alpha^2].
\end{equation}
From \eqref{SC} and \eqref{Ricci6} we obtain 

\begin{equation}\label{Coff3}
a_3 = \frac{r}{n-1} -(\xi \alpha) -\alpha^2, 
\end{equation}
and
\begin{equation}\label{Coff4}
b_3 = \frac{r}{n-1} -n[(\xi \alpha) +\alpha^2].
\end{equation}
Thus we have the following:
\begin{thm}
If an LTFSM of dimension \( n > 2\) is quasi-Einstein, then its associated scalars are given in \eqref{Coff3} and \eqref{Coff4}.
\end{thm}
We now suppose that an $LTFSM$ of dimension \( n > 2\) is generalized quasi-Einstein in the sense of Chaki \cite{Chaki01}. Then its Ricci tensor is of the form \eqref{Ricci4}. Hence for an $LTFSM$, \eqref{Ricci4} can be expressed as 
\begin{equation}\label{Ricci5}
S(Y,Z)=a g(Y,Z)+b \eta (Y)\eta (Z)+c[\zeta (Y)\eta (Z)+\zeta (Z)\eta (Y)],
\end{equation}
where $A(Y)=g(Y,\mu_1)=g(Y,\xi)=\eta(Y)$ and $\zeta(Y)=g(grad\alpha,Y)=(Y\alpha)$ for all $Y$. Since $g(grad\alpha,\xi)=(\xi\alpha)\neq 0,$ we now take projection of $grad \alpha$ onto the spatial distribution orthogonal to $\xi$ such that $grad^\perp \alpha=grad~\alpha+(\xi\alpha)\xi$ and hence $g(grad^\perp \alpha,\xi)=0$ and $||grad^\perp \alpha||^2=||grad~ \alpha||^2+(\xi\alpha)^2.$ We take $\mu_2=\frac{grad^\perp \alpha}{||grad^\perp \alpha||}.$ Then $g(\mu_2,\mu_2)=1$ and $\mu_2$ is orthogonal to $\mu_1=\xi.$ Consequently, \eqref{Ricci5} takes the form 
\begin{equation}\label{Ricci8}
S(Y,Z)=a g(Y,Z)+[b-2c(\xi\alpha)]\eta (Y)\eta (Z)=c||grad^\perp\alpha||[B (Y)\eta (Z)+B (Z)\eta (Y)],
\end{equation}
where $B(Y)=g(Y,~\mu_2)$ for all $Y.$ From \eqref{Ricci8} we get 
\begin{equation}\label{SC6}
r = n a - b+2c(\xi\alpha).
\end{equation}
 Also by using \eqref{Ricci} in \eqref{Ricci8} entails 
 \begin{equation}\label{Coeff8}
 (n-1)[(\xi\alpha)+\alpha^2]=a-b+2c(\xi\alpha).
 \end{equation}
 From \eqref{SC6} and \eqref{Coeff8} we obtain
 \begin{equation}\label{Coff7}
 a = \frac{r}{n-1} -(\xi \alpha) -\alpha^2, 
 \end{equation}
 and
 \begin{equation}\label{Coff8}
 b = \frac{r}{n-1} -n[(\xi \alpha) +\alpha^2]+2c(\xi\alpha).
 \end{equation}
 Again, setting $Z=\xi$ in \eqref{Ricci8} and then using \eqref{Ricci} we get
 \begin{equation}\label{Coef10}
 -(n-2)(Y\alpha)+c||grad^\perp\alpha||B(Y)=[a-b+(\xi\alpha)-(n-1)\alpha^2]\eta(Y).
 \end{equation}
 In view of \eqref{Coeff8} and \eqref{Coef10} yields (as $\xi\alpha\neq 0$)\\
 \begin{equation}\label{Coff13}
 c=1.
 \end{equation}
 Hence \eqref{Coff8} takes the form 
 \begin{equation}\label{Coff11}
  b = \frac{r}{n-1} -(n-2)(\xi \alpha) -n\alpha^2.
  \end{equation}
   This leads to the following:
\begin{thm}
If an LTFSM of dimension \( n > 2\) is generalized quasi-Einstein in the sense of Chaki, then its associated scalars are given in \eqref{Coff7}, \eqref{Coff13} and \eqref{Coff11}.
\end{thm}

\begin{cor}\cite{Shaikh03}
If an \(LCSM\) of dimension \(n > 2\) is quasi-Einstein or generalized quasi-Einstein (in the sense of Chaki), then its associated scalars are given by
$$a=\frac{r}{n-1}-(\rho-\alpha^2),$$
and
$$b=\frac{r}{n-1}-n(\rho-\alpha^2)$$
$\rho$ being the certain scalar function.
\end{cor}
\begin{cor}\cite{DSS2002}
Let $\left(M^n,\phi,\eta,\xi,g\right)$ be an $n$-dimensional (\(n > 2\)) LP-Sasakian manifold with coefficient $\alpha$. If $M^n$ is Einstein, then its scalar curvature is given by
\[
r = n\left\{(n-1)\alpha^2-\rho\Psi\right\},
\]
where $\rho=-\alpha(\xi)=-(\xi \alpha)$ and $\Psi=\operatorname{trace}(\phi)$.
\end{cor}
\begin{cor}\cite{DSS2002}
Let $(M^n,\phi,\eta,\xi,g)$ be an $n$-dimensional (\(n > 2\)) quasi-Einstein LP-Sasakian manifold with coefficient $\alpha$. Then the associated scalars $a$ and $b$ are completely determined by the scalar curvature $r$ and are given by
\[
a=\frac{r}{n-1}+\frac{\Psi \rho}{n-1}-\alpha^{2},
\]
and
\[
b=\frac{r}{n-1}-n\alpha^{2}+\frac{n\Psi \rho}{n-1},
\]
where $\rho=-\alpha(\xi)=-(\xi\alpha)$ and $\Psi=\operatorname{trace}(\phi)$.
\end{cor}
\section{Sectional Curvature of Lorentzian Torse-forming Structure Manifolds}
Let $(M^{n},g)$ be an $n$-dimensional Lorentzian manifold. A plane section at
a point $x\in M$ is a 2-dimensional vector subspace $E$ of the tangent space
$T_xM$. Such a plane section is said to be non-degenerate whenever the induced
metric on $E$ is non-degenerate. In particular, if $\{U,V\}$ constitutes a
basis of $E$, then
\[
g(U,U)g(V,V)-g(U,V)^2 \neq 0.
\]

The sectional curvature corresponding to the plane section $E$ at $x$ is given
by
\begin{equation}
K(x,E)
=
-\frac{g(R(U,V)U,V)}
       {g(U,U)g(V,V)-g(U,V)^2},
\label{eq:sectional_curvature}
\end{equation}
where $R$ denotes the Riemann curvature tensor of $(M,g)$. A non-degenerate
plane section is called timelike if it is generated by a timelike vector and a
spacelike vector.\\
\indent We now consider an $n$-dimensional $(n\geq3)$ $LTFSM$ with structure $(\phi,\eta,\xi,g)$ and suppose that $E$ is a $2$-dimensional linear subspace of the tangent space $T_xM$ at any point $x\in M$, and spanned by the unit timelike flow vector field $\xi$ and a spacelike vector $X$. Then since $E$ is timelike, the sectional curvature $K(x,E)$ of $E$ for any $X\in E$ is given by 
\begin{equation}
K(x,E)
=
\frac{g(R(\xi,X)\xi,X)}
       {g(X,X)+\eta(X)^2},
\label{eq:sectional_curvature1}
\end{equation}
which yield by virtue of \eqref{RC1} that $K(x,E)=\alpha^2+(\xi\alpha)$. This leads to the following:
\begin{thm}
Let \( M^{n}(\phi,\xi,\eta,g) \) be an \( n \)-dimensional (\(n \geq 3\)) $LTFSM$.  Then the sectional curvatue of each plane section spanned by the  timelike vector field \( \xi \) and a spacelike vector field is given by $$K(x,E)=\alpha^2+(\xi\alpha).$$ 
\end{thm}
\begin{cor}\cite{Shaikh09}
Let \( M^{n}(\phi,\xi,\eta,g) \) be an \( n \)-dimensional (\(n \geq 3\)) $LCSM$.  Then the sectional curvatue of each plane section spanned by the  timelike vector field \( \xi \) and a spacelike vector field is given by $$K(x,E)=\alpha^2-\rho,~\rho=-(\xi\alpha).$$ 

\end{cor}
\begin{cor}\cite{DSS2002}
Let \( M^{n}(\phi,\xi,\eta,g) \) be an \( n \)-dimensional (\(n \geq 3\)) LP-Sasakian manifold with a non-constant coefficient $\alpha$.  Then the sectional curvatue of each plane section spanned by the  timelike vector field \( \xi \) and a spacelike vector field is given by $$K(x,E)=\alpha^2-\rho,~\rho=-(\xi\alpha).$$ 

\end{cor}
\begin{cor}
In an $n$-dimensional (\(n \geq 3\)) LP-$\alpha$-Sasakian manifold $M^{n}(\phi,\xi,\eta,g)$, the sectional curvature corresponding to every plane section spanned by the timelike vector field $\xi$ and a spacelike vector field is given by
\[
K(x,E)=\alpha^{2}.
\]
\end{cor}
\begin{cor}
Let \( M^{n}(\phi,\xi,\eta,g) \) be an \( n \)-dimensional (\(n \geq 3\)) LP-Sasakian manifold with a constant coefficient $\alpha$.  Then the sectional curvatue of each plane section spanned by the  timelike vector field \( \xi \) and a spacelike vector field is given by $$K(x,E)=\alpha^2=constant.$$ 

\end{cor}
\begin{cor}\cite{Matsumoto89}
Let \( M^{n}(\phi,\xi,\eta,g) \) be an \( n \)-dimensional (\(n \geq 3\)) LP-Sasakian manifold.  Then the sectional curvatue of each plane section spanned by the  timelike vector field \( \xi \) and a spacelike vector field is given by $$K(x,E)=1.$$ 
\end{cor}
\begin{cor}
In an $n$-dimensional (\(n \geq 3\)) LP-$\alpha$-Kenmotsu manifold $M^{n}(\phi,\xi,\eta,g)$, the sectional curvature corresponding to every plane section spanned by the timelike vector field $\xi$ and a spacelike vector field is given by
\[
K(x,E)=-\alpha^{2}\neq \text{constant}.
\]
\end{cor}

\section{Some Examples of Lorentzian Torse-forming Structure Manifolds}
\begin{exm}
We consider the Friedmann–Lemaître–Robertson–Walker (FLRW) metric with coordinates $(t, r,\theta,\phi)$ and is given by
\[
ds^2 = -c^2 dt^2 + a^2(t)\left[ \frac{dr^2}{1 - k r^2} + r^2\bigl(d\theta^2 + \sin^2\!\theta\, d\phi^2\bigr) \right],
\]
where $c$ denotes the speed of light in a vacuum, $t$ represents cosmic time, and $r$ is the radial coordinate. The function $a(t)$, known as the scale factor, describes how the universe evolves over time by indicating whether it is expanding or contracting, and it is assumed to be a smooth, non-zero function of $t$. The spatial geometry of the universe is characterized by the parameter $k$: $k=+1$ corresponds to a closed universe with spherical geometry, $k=0$ describes a spatially flat and infinite universe, and $k=-1$ represents an open universe with hyperbolic geometry.\\
Let us consider a smooth vector field $\xi$ for this spacetime given as
\[
\xi=\left\{-\frac{1}{c},\,0,\,0,\,0\right\}.
\]

Then $\xi$ is unit timelike vector such that the associated covector field is obtained from the relation
$
\eta_j=g_{ij}\xi^i,
$
as
$
\eta=\{c,\,0,\,0,\,0\},
$
because
\[
g_{ij}\xi^i=-c^2\left(-\frac{1}{c}\right)=c.
\]

Hence, we obtain
\begin{equation}
\nabla_i\eta_j
=
-\frac{a'(t)}{c\,a(t)}
\left(
g_{ij}
+
\eta_i\eta_j
\right).
\end{equation}
where $a'(t)$ is derivative with respect to $t.$ Thus if we take 

$
\xi
=
-\frac{1}{c}\,
\frac{\partial}{\partial t}.
$
where $\frac{\partial}{\partial t}=
(-1,0,0,0)$, then from $(5.1)$ it follows that $\xi$ satisfies \eqref{UTT2}, where $\alpha=-\frac{a'(t)}{c\,a(t)}.$ This implies that $\xi$
is a unit timelike torse-forming vector field. Hence FLRW spacetime with $\xi
=
-\frac{1}{c}\,
\frac{\partial}{\partial t}$ is an example of $(LTFM)_4$.
\end{exm}
\begin{exm}\cite{MM16}
Generalized Robertson--Walker (GRW) spacetimes form a natural extension of Robertson--Walker spacetimes by relaxing the requirement of spatial homogeneity. An $n$-dimensional GRW manifold is a warped product Lorentzian manifold whose metric is given by
\[
ds^{2}=-dt^{2}+f^{2}(t)\,g^{*}_{\mu\nu}(x)\,dx^{\mu}dx^{\nu},
\]
where $f(t)$ is a positive function of time and $M^{*}$ is an arbitrary $(n-1)$-dimensional Riemannian manifold endowed with metric $g^{*}$. Unlike the Robertson--Walker case, no constant curvature condition is imposed on the spatial factor $M^{*}$.

The hypersurfaces $\{t\}\times M^{*}$ are totally umbilical, and the mean curvature of each slice is determined by
\[
H=\frac{\dot{f}}{f},
\]
which coincides with the Hubble parameter in cosmological models.

A fundamental characterization of GRW spacetime was established by Chen \cite{Chen2014}, who proved that a Lorentzian manifold is a GRW spacetime if and only if it admits a timelike concircular vector field $X$, satisfying
\[
X^{k}X_{k}<0,
\qquad
\nabla_{i}X_{j}=\nu\,g_{ij},
\]
for some scalar function $\nu$. The associated unit timelike vector field
\[
\eta_{k}=\frac{X_{k}}{\sqrt{-X^{2}}}
\]
is concircular and satisfies
\[
\nabla_{i}\eta_{j}
=\alpha\left(\eta_{i}\eta_{j}+g_{ij}\right),
\]
where $\alpha$ is a scalar function and $
\eta_j=g_{ij}\xi^i.
$ Hence a GRW manifold is an example of $LCSM.$
\end{exm}
\begin{exm}\cite{MM17}
The twisted cubic manifold represent the most general class of Lorentzian manifolds, extending the Robertson-Walker (RW) and generalized Robertson-Walker (GRW) spacetimes by allowing the scale factor \(f\) to depend on both time and position. In \cite{MM17} Mantica and Molinari established that the twisted cubic manifold admits a unit timelike torse-forming vector field, and also a Lorentzian manifold admitting a unit timelike torse-forming vector field is a twisted cubic manifold. A twisted cubic spacetime is a connected time oriented $4$-dimensional Lorentzian manifold \((M,g)\) equipped with the metric
\[
ds^{2} = -dt^{2} + f^{2}(t,\mathbf{x}) \left( dx^{2} + dy^{2} + dz^{2} \right),
\]
where \(f(t,\mathbf{x}) = 1 + \epsilon (t + x^2 + y^2 + z^2)\) with \(0 < \epsilon \ll 1\). This metric represents a small inhomogeneous perturbation of the flat Robertson-Walker metric, with the scale factor depending on both time and spatial coordinates, thus providing a concrete example of a twisted spacetime.

For this twisted spacetime, we define the vector field \(u\) by \(u_0 = -1\), \(u_{\mu} = 0\), which satisfies \(u^{i}u_{i} = -1\), making it a unit timelike vector field. The Christoffel symbols for this metric, as derived in the appendix of Mantica and Molinari \cite{MM17}, are given by \(\Gamma_{00}^{k} = 0\), \(\Gamma_{\mu 0}^{\nu} = (\dot{f}/f)\delta_{\mu}^{\nu}\), \(\Gamma_{\mu\nu}^{0} = f\dot{f}g_{\mu\nu}^{*}\), and \(\Gamma_{\mu\nu}^{\rho} = \Gamma_{\mu\nu}^{*\rho} + (f_{\nu}/f)\delta_{\mu}^{\rho} + (f_{\mu}/f)\delta_{\nu}^{\rho} - (f^{\rho}/f)g_{\mu\nu}^{*}\), where \(\dot{f} = \partial_t f\), \(f_{\mu} = \partial_{\mu}f\), and \(f^{\mu} = g^{*\mu\nu}f_{\nu}\). Computing the covariant derivative of \(u\), we find \(\nabla_{0}u_{0} = 0 = \alpha(u_{0}^{2} + g_{00})\), \(\nabla_{\mu}u_{0} = 0\), \(\nabla_{0}u_{\mu} = 0\), and \(\nabla_{\mu}u_{\nu} = f\dot{f}\delta_{\mu\nu}\). This yields \(\nabla_{i}u_{j} = \alpha(u_{i}u_{j} + g_{ij})\) with \(\alpha = \dot{f}/f = \epsilon/[1 + \epsilon(t + x^2 + y^2 + z^2)]\), confirming that \(u\) is a unit timelike torse-forming vector field.
The uniqueness of this torse-forming vector field is established through a rigorous argument. Suppose there exists another unit timelike torse-forming vector field \(v\) satisfying \(\nabla_{i}v_{j} = \lambda(v_{i}v_{j} + g_{ij})\). In the twisted frame, we write \(v_0 = \cosh a\), \(v_{\mu} = V_{\mu}\sinh a\), with \(a(t,\mathbf{x}) \neq 0\) and \(g^{*\mu\nu}V_{\mu}V_{\nu} = f^{2}\). The equations \(\nabla_{0}v_{0} = \lambda(v_{0}^{2} - 1)\) and \(\nabla_{0}v_{\mu} = \lambda v_{0}v_{\mu}\) yield \(\dot{a}\sinh a = \lambda\sinh^{2}a\) and \(\sinh a(f\dot{V}_{\mu} + V_{\mu}\dot{f}) = 0\). If \(a \neq 0\), then \(fV_{\mu}\) is independent of \(t\), which is compatible with \(g^{*\mu\nu}V_{\mu}V_{\nu} = f^{2}\) only if \(\dot{f} = 0\). Thus, for \(\dot{f} \neq 0\), we must have \(a = 0\), implying \(v = u\) up to sign. Therefore, the torse-forming vector field is unique in non-trivial twisted spacetimes.

The geometric significance of this torse-forming vector field lies in its role as the comoving observer field. The function \(\alpha = \dot{f}/f\) represents the mean curvature of the spatial slices, which in a twisted spacetime is not uniform: \(H(t,\mathbf{x}) = \alpha(t,\mathbf{x}) = u^{k}\nabla_{k}\log f\). This non-uniform expansion gives rise to an energy flux and anisotropic stress in the corresponding Einstein field equations.
Hence twisted cubic spacetime is an $(LTFM)_4.$
\end{exm}
\begin{exm}\cite{Shaikh09}
Consider the 3-dimensional manifold $M=\{(x,y,z)\in\mathbb{R}^3\},$ where $(x,y,z)$ denotes the standard coordinates on $\mathbb{R}^3$. Let $\{e_1,e_2,e_3\}$ be the linearly independent global frame on $M$ such that \[ e_1=e^{-z}\left(\frac{\partial}{\partial x} +y\frac{\partial}{\partial y}\right),\qquad e_2=e^{-z}\frac{\partial}{\partial y}, e_3=e^{-2z}\frac{\partial}{\partial z}. \] We define a Lorentzian metric $g$ on $M$ by $ g(e_1,e_1)=g(e_2,e_2)=1, g(e_3,e_3)=-1, $ $g(e_1,e_2)=g(e_1,e_3)=g(e_2,e_3)=0.$
Thus in terms of coordinates, the metric can be written as $$g = e^{2z}(1 + y^2) dx^2 + e^{2z} dy^2 - e^{4z} dz^2 - 2ye^{2z} dx dy.$$ Thus, $\{e_1,e_2,e_3\}$ forms a pseudo-orthonormal frame, with $e_3$ being the timelike vector. Let $\eta$ be the $1$-form associated with $e_3$, defined by $ \eta(U)=g(U,e_3), U\in\chi(M),$  and let $\phi$ be the $(1,1)$-tensor field defined by  $\phi e_1=e_1,\qquad \phi e_2=e_2,\qquad \phi e_3=0.$ In terms of local coordinates, the matrix representation of $\phi$ is given as 
$[\phi] = \begin{pmatrix} 1 & 0 & 0 \\ 0 & 1 & 0 \\ 0 & 0 & 0 \end{pmatrix}$.

 By extending these relations linearly, we obtain, for arbitrary vector fields $U,W\in\chi(M)$, $ \eta(e_3)=-1,  \phi^2U=U+\eta(U)e_3=\phi U,~ g(\phi U,\phi W) =g(U,W)+\eta(U)\eta(W).$ Hence, by taking $\xi=e_3$, the structure $(\phi,\xi,\eta,g)$ satisfies the defining relations of a Lorentzian paracontact structure on $M$. The Lie brackets of the frame fields are readily computed as \[ [e_1,e_2]=-e^{-z}e_2,\qquad [e_1,e_3]=e^{-2z}e_1,\qquad [e_2,e_3]=e^{-2z}e_2. \] Using these relations together with the Koszul formula for the Lorentzian metric $g$, the Levi-Civita connection $\nabla$ is given by \begin{align*} \nabla_{e_1}e_3 &= e^{-2z}e_1, & \nabla_{e_1}e_2 &= 0, & \nabla_{e_1}e_1 &= e^{-2z}e_3,\\ \nabla_{e_2}e_3 &= e^{-2z}e_2, & \nabla_{e_2}e_2 &= e^{-2z}e_3-e^{-z}e_1, & \nabla_{e_2}e_1 &= e^{-2z}e_2, \end{align*}  \[ \nabla_{e_3}e_3=0,\qquad \nabla_{e_3}e_2=0,\qquad \nabla_{e_3}e_1=0. \] These connection formulas show that the structure $(\phi,\xi,\eta,g)$ realizes that $\xi=e_3$ is unit torse-forming timelike vector field such that \eqref{UTT2} and \eqref{Omega} holds for  $\alpha=e^{-2z}, $ which is nonzero everywhere on $M$, and satisfies  $(X\alpha)=\rho\,\eta(X), $ $\mu=2e^{-4z}.$
Therefore $\xi$ is unit timelike torse-forming as well as concircular vector field and hence $M^3(\phi,\xi,\eta,g)$ is an example of $(LTFSM)_3$ as well as $(LCSM)_3.$
\end{exm}
\begin{exm}\cite{Shaikh09}
Consider the $4$-dimensional manifold $M=\left\{(x,y,z,t)\in\mathbb{R}^4\mid t\neq 0\right\},$ where $(x,y,z,t)$ are the standard coordinates on $\mathbb{R}^4$. Let $\{e_1,e_2,e_3,e_4\}$ be a linearly independent global frame on $M$ defined by \[ e_1=t\left(\frac{\partial}{\partial x} +y\frac{\partial}{\partial y}\right),\qquad e_2=t\frac{\partial}{\partial y}, \] \[ e_3=t\left(\frac{\partial}{\partial y} +\frac{\partial}{\partial z}\right),\qquad e_4=t^3\frac{\partial}{\partial t}. \] Define the $(1,1)$-tensor field $\phi$, the vector field $\xi$, the $1$-form $\eta$, and the Lorentzian metric $g$ by \[ \phi e_1=e_1,\qquad \phi e_2=e_2,\qquad \phi e_3=e_3,\qquad \phi e_4=0, \] and \[ \xi=(t)^3\frac{\partial}{\partial t},\qquad \eta(X)=g(X,e_4), \phi X=X+\eta(X)\xi=\phi^2 X \qquad X\in\chi(M). \] 
In terms of local coordinates, the matrix representation of $\phi$ is given as 
$$[\phi] = \begin{pmatrix} 
1 & 0 & 0 & 0 \\ 
0 & 1 & 0 & 0 \\ 
0 & 0 & 1 & 0 \\ 
0 & 0 & 0 & 0 
\end{pmatrix}.$$
The metric $g$ is determined on the chosen frame by \[ g(e_1,e_1)=g(e_2,e_2)=g(e_3,e_3)=1, g(e_4,e_4)=-1, \]  $ g(e_i,e_j)=0, i\neq j, i,j=1,2,3,4.$ In local coordinates, $g$ can be written as
 $$g = \frac{1}{t^2} \left[ (1+y^2)dx^2 + dy^2 + 2dz^2 - 2y\,dx\,dy + 2y\,dx\,dz - 2\,dy\,dz \right] - \frac{1}{t^6} dt^2.$$
  Thus, $\{e_1,e_2,e_3,e_4\}$ forms a pseudo-orthonormal frame, with $e_4$ being timelike. The non-vanishing Lie brackets of the frame fields are given by \[ [e_1,e_2]=-te_2,\qquad [e_1,e_4]=-(t)^2e_1, \] \[ [e_2,e_4]=-(t)^2e_2,\qquad [e_3,e_4]=-(t)^2e_3. \] Identifying $e_4$ with $\xi$ and applying the Koszul formula to the Lorentzian metric $g$, the corresponding Levi--Civita connection $\nabla$ is obtained as follows: \begin{align*} \nabla_{e_1}e_4 &= -(t)^2e_1, & \nabla_{e_2}e_1 &= te_2, & \nabla_{e_1}e_1 &= -(t)^2e_4,\\ \nabla_{e_2}e_4 &= -(t)^2e_2, & \nabla_{e_3}e_4 &= -(t)^2e_3, & \nabla_{e_3}e_3 &= -(t)^2e_4, \end{align*}  $ \nabla_{e_2}e_2=-(t)^2e_4-te_1.$ The above connection relations show that the structure $(\phi,\xi,\eta,g)$ satisfies the defining conditions of an $(LTFSM)_4$ as well as $(LCSM)_4$, where the associated function $\alpha$ is given by  $\alpha=-(t)^2\neq0,$ and it satisfies $ (X\alpha)=\mu\,\eta(X), $  $ \mu=2(t)^4. $ Therefore, the vector field $\xi$ is a unit timelike torse-forming vector field and, in particular, is concircular. Consequently, $ M^4(\phi,\xi,\eta,g) $ is an example of $(LTFSM)_4$ and $(LCSM)_4$. 
\end{exm}
\begin{exm}
Consider the 3-dimensional manifold $M = \{(x,y,z) \in \mathbb{R}^3\}$, where $(x,y,z)$ denote the standard coordinates on $\mathbb{R}^3$. Let $\{e_1, e_2, e_3\}$ be a global frame on $M$ defined by
\[
e_1 = e^{-x}\frac{\partial}{\partial x}, \quad 
e_2 = \frac{\partial}{\partial y}, \quad 
e_3 = \frac{\partial}{\partial z}
,\]
and let $\phi$ be the $(1,1)$-tensor field defined by \[ \phi e_1=e_1,\qquad \phi e_2=e_2,\qquad \phi e_3=0. \]
In terms of local coordinates, the matrix representation of $\phi$ is given as 
$$[\phi] = \begin{pmatrix} 1 & 0 & 0 \\ 0 & 1 & 0 \\ 0 & 0 & 0 \end{pmatrix}.$$
Define the Lorentzian metric $g$ on $M$ such that
\[
g(e_1,e_1) = 1, \quad g(e_2,e_2) = 1, \quad g(e_3,e_3) = -1, \quad g(e_i, e_j) = 0 \quad (i \neq j).
\]
Set $\xi = e_3$ as the characteristic vector field, with dual 1-form\\ $\eta(X)=g(X,e_4), \phi X=X+\eta(X)\xi=\phi^2 X,~ X\in\chi(M).$

By using Koszul formula to the Lorentzian metric $g$, the corresponding Levi--Civita connection $\nabla$ is obtained as follows:
\[
\nabla_{e_1} e_3 = e^{-x} e_1, \quad \nabla_{e_2} e_3 = e^{-x} e_2, \quad \nabla_{e_3} e_3 = 0.
\]
This implies that the flow vector field $\xi$ satisfies relation:
$
\nabla_X \xi = \alpha \{X + \eta(X)\xi\},
$
where $\alpha = e^{-x} \neq 0$.

    Since $g(\xi, \xi) = -1$ and $\nabla_X \xi = \alpha X + \tau(X)\xi$ with $\tau(X) = \alpha \eta(X)$, $\xi$ is a unit timelike torse-forming vector field.\\
     For $\xi$ to be concircular, the condition $X\alpha = \mu \eta(X)$ must hold for all $X$. Evaluating at $e_1$ we get
    \[
    e_1(\alpha)=(e_1\alpha)=\nabla_{e_1}\alpha = \left(e^{-x} \frac{\partial}{\partial x}\right)(e^{-x}) = -e^{-2x} \neq 0.
    \]
    However, $\mu \eta(e_1) = \mu g(e_1, e_3) = 0$. Since $e_1(\alpha) \neq \mu \eta(e_1)$, the condition fails. 
Therefore, $\xi$ is unit timelike torse-forming but not concircular. Hence $ M^3(\phi,\xi,\eta,g) $ is an example of $(LTFSM)_3$ but not $(LCSM)_3$.
\end{exm}
\begin{exm}
Consider the 4-dimensional manifold $M = \{(x, y, z, t) \in \mathbb{R}^4 \mid t \neq 0\}$. where $(x,y,z,t)$ denote the standard coordinates on $\mathbb{R}^4$. Let $\{e_1, e_2, e_3, e_4\}$ be a global frame given by
\[
e_1 = e^{-x}\left(\frac{\partial}{\partial x} + y\frac{\partial}{\partial y}\right), \quad 
e_2 = e^{-x}\frac{\partial}{\partial y}, \quad 
e_3 = e^{-x}\left(\frac{\partial}{\partial y} + \frac{\partial}{\partial z}\right), \quad 
e_4 = e^{-x}\frac{\partial}{\partial t}.
\]

Define the Lorentzian metric $g$ on $M$ such that
\[
g(e_1, e_1) = g(e_2, e_2) = g(e_3, e_3) = 1, \quad g(e_4, e_4) = -1, \quad g(e_i, e_j) = 0 \quad (i \neq j).
\]
In coordinates $g$ can be written as
$$g = e^{2x} \left[ (1+y^2)dx^2 + dy^2 + 2dz^2 - 2y\,dx\,dy + 2y\,dx\,dz - 2\,dy\,dz \right] - e^{2x} dt^2,$$
and let $\phi$ be the $(1,1)$-tensor field defined by \[ \phi e_1=e_1,\qquad \phi e_2=e_2,\qquad \phi e_3=e_3,\qquad \phi e_4=0. \]
The matrix representation of $\phi$ is given as 
$$[\phi] = \begin{pmatrix} 
1 & 0 & 0 & 0 \\ 
0 & 1 & 0 & 0 \\ 
0 & 0 & 1 & 0 \\ 
0 & 0 & 0 & 0 
\end{pmatrix}.$$
Set $\xi = e_4$, $\xi = e_4$ as the characteristic vector field, with dual 1-form $\eta(X)=g(X,e_4), \phi X=X+\eta(X)\xi=\phi^2 X,~ X\in\chi(M),$ by using Koszul formula to the Lorentzian metric $g$, the corresponding Levi--Civita connection $\nabla$ is obtained as follows:
\[
\nabla_{e_1} e_4 = -e^{-x} e_1, \quad \nabla_{e_2} e_4 = -e^{-x} e_2, \quad \nabla_{e_3} e_4 = -e^{-x} e_3, \quad \nabla_{e_4} e_4 = 0
\]

The covariant derivative takes the form $\nabla_X \xi = \alpha \{X + \eta(X)\xi\}$ with $\alpha = -e^{-x} \neq 0$.

Now, $g(\xi, \xi) = -1$ and $\nabla_X \xi = \alpha X + \alpha \eta(X)\xi$, $\xi$ is a unit timelike torse-forming vector field. Now
    \[
    e_1(\alpha)=(e_1\alpha) = \left(e^{-x} \frac{\partial}{\partial x}\right)(-e^{-x}) = e^{-2x} \neq 0.
    \]
    Since $e_1$ is orthogonal to $\xi$, $\eta(e_1) = 0$, giving $\mu \eta(e_1) = 0$. The condition $(X\alpha) = \mu \eta(X)$ fails because $e_1(\alpha) \neq 0$.
Thus, $\xi$ is unit timelike torse-forming but not concircular. Hence $ M^4(\phi,\xi,\eta,g) $ is an example of $(LTFSM)_4$ but not $(LCSM)_4$.
\end{exm}
\begin{exm}
Consider the $5$-dimensional manifold $M^5 = \{(x,y,w,z,t) \in \mathbb{R}^5 \mid t \neq 0\}$, where $(x,y,w,z,t)$ denote the standard coordinates on $\mathbb{R}^5$. Let $f(x)$ be a smooth, non-constant function depending only on $x$ (for example, $f(x) = e^{-x}$). We define a linearly independent global frame $\{e_1, e_2, e_3, e_4, e_5\}$ on $M$ by
\[
e_{1}=e^{-f(x)}\left(\frac{\partial }{\partial x}+y\frac{\partial }{\partial y}\right),
e_{2}=e^{-f(x)}\frac{\partial }{\partial y},
e_{3}=e^{-f(x)}\frac{\partial }{\partial w},
e_{4}=e^{-f(x)}\frac{\partial }{\partial z}, 
e_{5}=e^{-f(x)}\frac{\partial }{\partial t}=:\xi.
\]
Let $g$ be the standard pseudo-orthonormal Lorentzian metric defined on this frame given by
\[
g(e_{i},e_{i})=1 \quad (i=1,2,3,4),\qquad g(e_{5},e_{5})=-1,\qquad g(e_{i},e_{j})=0 \quad (i\ne j).
\]
In terms of local coordinates the metric can be expressed as
$$g = e^{2f(x)} \left[ (1+y^2)dx^2 + dy^2 + dw^2 + dz^2 - 2y\,dx\,dy \right] - e^{2f(x)} dt^2.$$
In matrix form $(g_{ij})$ with respect to the standard coordinate basis $\{dx, dy, dw, dz, dt\}$, the metric tensor is expressed as:
\[
g_{ij} = \begin{pmatrix} 
e^{2f(x)}(1+y^2) & -y e^{2f(x)} & 0 & 0 & 0 \\ 
-y e^{2x} & e^{2f(x)} & 0 & 0 & 0 \\ 
0 & 0 & e^{2f(x)} & 0 & 0 \\ 
0 & 0 & 0 & e^{2f(x)} & 0 \\ 
0 & 0 & 0 & 0 & -e^{2f(x)} 
\end{pmatrix}.
\]
Thus, $\{e_1, e_2, e_3, e_4, e_5\}$ forms a pseudo-orthonormal frame, with $\xi = e_5$ being the unit timelike vector field. Let $\eta$ be the $1$-form associated with $\xi$, defined by $\eta(U)=g(U,\xi)$ for any $U \in \chi(M)$.
Hence
$
g(\xi,\xi)=g(e_5,e_5)=-1,
$
so $\xi$ is a unit timelike vector field. Define the dual 1-form $\eta$ by
$
\eta(X)=g(X,\xi),
$
and the $(1,1)$-tensor field $\phi$ by
\[
\phi e_i = e_i \quad (i=1,2,3,4),\qquad \phi e_5=0.
\]
The matrix representation of $\phi$ is given as 
$$[\phi] = \begin{pmatrix} 
1 & 0 & 0 & 0 & 0 \\ 
0 & 1 & 0 & 0 & 0 \\ 
0 & 0 & 1 & 0 & 0 \\ 
0 & 0 & 0 & 1 & 0 \\ 
0 & 0 & 0 & 0 & 0 
\end{pmatrix}.$$
Then $g(X,\xi)=\eta(X),\phi^2 = I + \eta\otimes \xi=\phi$, $\eta(\xi)=-1$, and $g(\phi X,\phi Y)=g(X,Y)+\eta(X)\eta(Y)$.
By applying the Koszul formula, the covariant derivatives with respect to the timelike vector field $\xi = e_5$ are given by\\
$
\nabla _{e_{1}}\xi =-f^{\prime }(x)e^{-f(x)}e_{1}, \quad 
\nabla _{e_{2}}\xi =-f^{\prime }(x)e^{-f(x)}e_{2}, \quad
\nabla _{e_{3}}\xi =-f^{\prime }(x)e^{-f(x)}e_{3},\\
\nabla _{e_{4}}\xi =-f^{\prime }(x)e^{-f(x)}e_{4}, \quad
\nabla _{e_{5}}\xi =0.\\
$
Hence it follows that for any vector field $X \in \chi(M)$, we get
\[
\nabla _{X}\xi =\alpha \{X+\eta (X)\xi \},
\]
where
$
\alpha =-f^{\prime }(x)e^{-f(x)}\ne 0
$
is a smooth function. Hence $\xi$ is unit timelike torse-forming vector field.

For $\xi$ to be concircular, the scalar function $\alpha$ must obey the strict differential condition
\[
X(\alpha )=\mu \,\eta (X),
\]
for some scalar function $\mu$. If $X = e_1$, then 
    \[
    e_{1}(\alpha )=\left(e^{-f(x)}\frac{\partial }{\partial x}\right)\left(-f^{\prime }(x)e^{-f(x)}\right)=e^{-f(x)}\left(-f^{\prime \prime }(x)e^{-f(x)}+[f^{\prime }(x)]^{2}e^{-f(x)}\right)\ne 0.
    \]
   Also Since $e_{1}$ is orthogonal to $\xi = e_5$,
    \[
    \eta (e_{1})=g(e_{1},e_{5})=0 \implies \mu \,\eta (e_{1})=0.
    \]
Thus $\xi$ is not concircular. Hence $ M^5(\phi,\xi,\eta,g) $ is an example of $(LTFSM)_5$ but not $(LCSM)_5$.
\end{exm}

%
\section{Conclusion}
This paper mainly introduces a large class of Lorentzian manifolds admitting a unit timelike torse-forming vector fields, called \emph{Lorentzian torse-forming structure manifold (briefly, LTFSM)} in which the structure tensor $\phi$ is aligned with the torse-forming vector field.\\
\indent In an $LTFSM$, the projection behaviour of $\phi$ acts as the geometric constraint that enforces spatial isotropy, filtering out shear and vorticity. Consequently, the torse-forming scalar $\alpha$ is forced to act exactly as the Hubble expansion parameter, mapping the manifolds' purely geometric Ricci curvature directly to the gravitational deceleration caused by matter density and pressure.\\

 \indent This class namely, $LTFSM$ contains  $LCSMs$, $LSP$-Sasakian manifolds, and special type of $LP$-Sasakian manifolds with coefficient $\alpha$, $LP$-$\alpha$-Sasakian manifolds, $LP$-Sasakian manifolds and $LP$-$\alpha$-Kenmotsu manifolds. In Fig 1 and Fig 2 we have explored the interrelations among such various classes of manifolds. It is interesting to note that in an $LCSM$, the Ricci operator $Q$ commutes with the structure tensor $\phi$ \cite{Shaikh09}, but it is not true in an $LTFSM$. However, if the gradient of the of associated scalar $\alpha$ is codirectional with the torse-forming vector field $\xi,$ then the Ricci operator of $LTFSM$ commutes with the structure tensor $\phi$ and hence it reduces to an $LCSM.$ It is shown that a conformally flat \emph{LTFSM} is a generalized quasi-Einstein manifold in the sense of Chaki \cite{Chaki01}, and hence it is a manifold of generalized quasi-constant curvature, although a conformally flat $LCSM$ is a manifold of quasi-constant curvature. We determine the sectional curvature of an $LTFSM$ and observed that it is non-constant throughout the manifold. The sectional curvature of all timelike planes containing $\xi$ depends entirely on the behaviour of the scalar $\alpha$ along the direction of $\xi.$ It is completely independent of the choice of orthogonal vector $X$. If $\alpha$ is the along leaves of the distribution then these planes have a sectional curvature equal to $\alpha^2$.  Also we have deduced the associated scalars of quasi-Einstein and generalised quasi-Einstein $LTFSM$. Finally, we have exhibited several non-trivial examples which are only $LTFM$, $LTFSM$ as well as $LCSM$ and also $LTFSM$ but not $LCSM$.
\section{Future scope of research}
There is a wide scope of future research on $LTFSM$ such as to investigate different types of symmetry, semisymmety and pseudosymmetry, different types of Ricci solitons, curvature inheritance symmetry, various types of submanifolds etc.\\ 
 
\section*{Declarations}

 \subsection*{Data Availability Statement:}
 The manuscript has no associated data.

 \subsection*{Conflict of Interest}
 The author declare that he has no conflict of interest.

\textbf{Acknowledgement:}
 All the algebraic computations are performed by a program in Wolfram Mathematica developed by A. A. Shaikh. 
%

\end{document}